\input amstex
\documentstyle{amsppt}

\TagsOnRight \NoBlackBoxes \NoRunningHeads

\define\R{\Bbb R}
\define\Z{\Bbb Z}
\define\C{\Bbb C}

\define\X{\frak X}
\define\al{\alpha}
\define\be{\beta}

\define\la{\lambda}
\define\Om{\Omega}
\define\om{\omega}

\define\inr{\operatorname{in}}
\define\out{\operatorname{out}}

\define\wt{\widetilde}
\define\const{\operatorname{const}}
\define\Ga{\Gamma}

\define\Conf{\operatorname{Conf}}
\define\tht{\thetag}

\define\N{\Cal N}

\define\de{\delta}

\define\wh{\widehat}

\define\diag{\operatorname{diag}}
\define\SGN{\operatorname{SGN}}
\define\De{\Delta}
\define\Prob{\operatorname{Prob}}
\define\hypgeom{{\operatorname{hypergeom}}}
\define\askey{{\operatorname{Askey-Lesky}}}

\topmatter

\title Representation theory and random point
processes
\endtitle

\author Alexei Borodin and Grigori Olshanski
\endauthor

\abstract On a particular example we describe how to state and to
solve the problem of harmonic analysis for groups with
infinite--dimensional dual space. The representation theory for
such groups differs in many respects from the conventional theory.
We emphasize a remarkable connection with random point processes
that arise in random matrix theory. The paper is an extended
version of the second author's talk at the Congress.
\endabstract

\endtopmatter

\document

\head Introduction \endhead

In this paper we would like to discuss a connection between two
areas of mathematics which until recently seemed to be rather
distant from each other: (1) noncommutative harmonic analysis on
groups and (2) some topics in probability theory related to
random point processes. In order to make the paper accessible to
readers not familiar with either of these areas, we will explain
all needed basic concepts.

The purpose of harmonic analysis is to decompose natural
representations of a given group on irreducible representations.
By natural representations we mean those representations that
are produced, in a natural way, from the group itself. Examples
include the regular representation, which is realized in the
$L^2$ space on the group, or a quasiregular representation,
which is built from the action of the group on a homogeneous
space.

In practice, a natural representation often comes together with
a distinguished cyclic vector. Then the decomposition into
irreducibles is governed by a measure, which may be called the
{\it spectral measure.\/} The spectral measure lives on the dual
space to the group, the points of the dual being the irreducible
unitary representations. There is a useful analogy in analysis:
expanding a given function on eigenfunctions of a self--adjoint
operator. Here the spectrum of the operator is a counterpart of
the dual space.

If our distinguished vector lies in the Hilbert space of the
representation, then the spectral measure has finite mass and
can be normalized to be a probability measure. \footnote{It may
well happen that the distinguished vector belongs to an
extension of the Hilbert space (just as in analysis, one may
well be interested in expanding a function which is not square
integrable). For instance, in the case of the regular
representation of a Lie group one usually takes the delta
function at the unity of the group, which is not an element of
$L^2$. In such a situation the spectral measure is infinite.
However, we shall deal with finite spectral measures only.}

Now let us turn to {\it random point processes\/} (or random
point fields), which form a special class of stochastic
processes. In general, a stochastic process is a family of
random variables, while a point process (or random point field)
is a random point configuration. By a (nonrandom) point
configuration we mean an unordered collection of points in a
locally compact space $\frak X$. This collection may be finite
or countably infinite, but it cannot have accumulation points in
$\frak X$. To define a point process on $\frak X$, we have to
specify a probability measure on $\Conf(\frak X)$, the set of
all point configurations.

One classical example is the Poisson process, which is employed
in a lot of probabilistic models and constructions. Another
important example (or rather a class of examples) comes from
random matrix theory. Given a probability measure on a space of
$N\times N$ matrices, we pass to the matrix eigenvalues and thus
obtain a random $N$--point configuration. In a suitable scaling
limit transition (as $N\to\infty$), it turns into a point
process living on infinite point configurations.

As long as we are dealing with ``conventional'' groups (finite
groups, compact groups, real or $p$--adic reductive groups,
etc.), representation theory seems to have nothing in common
with point processes. However, the situation drastically changes
when we turn to ``big'' groups whose irreducible representations
depend on infinitely many parameters. Two basic examples are the
infinite symmetric group $S(\infty)$ and the
infinite--dimensional unitary group $U(\infty)$, which are
defined as unions of ascending chains of finite or compact
groups
$$
S(1)\subset S(2)\subset S(3)\subset\dots, \qquad U(1)\subset
U(2)\subset U(3)\subset\dots,
$$
respectively. It turns out that for such groups, the clue to the
problem of harmonic analysis can be found in the theory of point
processes.

The idea is to convert any infinite collection of parameters,
which corresponds to an irreducible representation, to a point
configuration. Then the spectral measure defines a point
process, and one may try to describe this process (hence the
initial measure) using appropriate probabilistic tools.

This approach was first applied to the group $S(\infty)$ (see
the surveys Borodin--Olshanski \cite{BO2}, Olshanski
\cite{Ol6}). In the present paper we discuss the group
$U(\infty)$, our exposition is mainly based on Olshanski
\cite{Ol7} and Borodin--Olshanski \cite{BO6}. Notice that the
point processes arising from the spectral measures do not
resemble the Poisson process but are close to the processes of
random matrix theory.

\subhead Acknowledgement\endsubhead This research was partially
conducted during the period the first author (A.~B.) served as a
Clay Mathematics Institute Research Fellow. He was also
partially supported by the NSF grant DMS-0402047. The second
author (G.~O.) was supported by the CRDF grant RM1-2543-MO-03.

\head 1. Dual space and the problem of harmonic analysis
\endhead

Recall that a {\it unitary representation\/} of a group $G$ in a
Hilbert space $H$ is a homomorphism of $G$ into the group of
unitary operators in $H$. For instance, if $G$ is a locally
compact topological group then there is a natural representation
generated by the (say, right) action of $G$ on itself, called
the {\it regular representation\/}. Its space is the $L^2$ space
formed with respect to the Haar measure on $G$, and the
operators of the representation are given by
$$
(R(g)f)(x)=f(xg), \qquad g\in G, \quad x\in G, \quad f\in
L^2(G). \tag1.1
$$

A unitary representation is said to be {\it irreducible\/} if it
is not a direct sum of other representations. Irreducible
representations are elementary objects like simple modules. A
general unitary representation $T$ is, in a certain sense, built
from irreducible ones: in simplest cases $T$ is decomposed into
a direct sum of irreducibles, and in more sophisticated
situations, direct sum is replaced by ``direct integral''.
\footnote{This claim is true under certain assumptions on the
group $G$ or on the representation $T$, but we don't want to
discuss technicalities here. Under additional (but still rather
broad assumptions), the decomposition into irreducibles is
essentially unique.}

Two fundamental problems of unitary representation theory are:

1. Given a group $G$, find all its irreducible unitary
representations.

2. For most natural representations of $G$ (e.g., the regular
representation), describe their decomposition on irreducibles.

The set of (equivalence classes of) irreducible unitary
representations of $G$ is called the {\it dual space\/} to $G$
and is denoted by $\wh G$. Thus, the first problem is the
description of $\wh G$. The second problem is called the {\it
problem of harmonic analysis\/}. It can be viewed as a
noncommutative generalization of the classical Fourier analysis.

These two problems were extensively studied for ``conventional''
groups. The existing literature is immense, and surveying it is
beyond the scope of the present paper. What is important for us
is that both problems, with appropriate refinement, make sense
for certain ``nonconventional'' groups as well. These are the
groups of automorphisms of infinite--dimensional Riemannian
symmetric spaces and also certain combinatorial analogs of such
groups, which are built with the help of the infinite symmetric
group.

Results on construction and classification of irreducible
representations for the automorphism groups and their
combinatorial analogs can be found in Olshanski \cite{Ol1},
\cite{Ol5}, \cite{Ol2}, \cite{Ol4}, Pickrell \cite{Pi2},
Nessonov \cite{Nes}. The construction of natural reducible
representations for these groups and related questions are
discussed in Pickrell \cite{Pi1}, Kerov--Olshanski--Vershik
\cite{KOV1}, \cite{KOV2}, Olshanski \cite{Ol7}. In the present
paper we focus on a single group $G$, which is $U(\infty)\times
U(\infty)$. The reason why we consider not the group $U(\infty)$
but the product of its two copies will be explained below. Here
we would only like to note that $U(\infty)$ (or an appropriate
completion thereof) can be viewed as an infinite--dimensional
Riemannian symmetric space, and then $U(\infty)\times U(\infty)$
arises as a group of automorphisms of that space.

\head 2. The dual space $\wh{U(N)}$ and spherical
representations of $U(N)\times U(N)$
\endhead

In this section we briefly describe a few necessary facts about
representations of the groups $U(N)$. The material is classical,
\footnote{See, e.g., Weyl \cite{We}, Zhelobenko \cite{Zhe},
Helgason \cite{He}.} we present it in a form which will help to
understand the subsequent infinite--dimensional generalization.

For $N=1,2,\dots$ let $U(N)$ denote the group of unitary
matrices of size $N\times N$. This group is compact. Its
irreducible representations are parametrized by {\it signatures
of length $N$\/}, that is, $N$--tuples $\la=(\la_1,\dots,\la_N)$
of integers such that $\la_1\ge\dots\ge\la_N$. \footnote{Another
term for collections $\la$ is ``dominant highest weights for
$U(N)$''.} Thus, the dual space $\wh{U(N)}$ can be viewed as a
countable discrete subset of $\R^N$.

Let $\pi^\la$ denote the irreducible representation
corresponding to a signature $\la\in\wh{U(N)}$, $\dim\pi^\la$
denote the dimension of the representation space, and $R_N$ be
the regular representation of $U(N)$ in the Hilbert space
$L^2(U(N))$. The decomposition of $R_N$ looks as follows
$$
R_N=\bigoplus_{\la\in\wh{U(N)}} \dim\pi^\la\cdot\pi^\la
$$
In other words, each irreducible representation enters the
regular representation with multiplicity equal to the dimension
of this irreducible representation. This is a special case of a
general result valid for any compact group, the Peter--Weyl
theorem.

We observe now that the group $U(N)$ acts on itself both on the
right and on the left, so that $U(N)$ becomes a homogeneous
space $U(N)\times U(N)/\diag(U(N))$, where $\diag(U(N))$ stands
for the diagonal subgroup in $U(N)\times U(N)$. This enables us
to extend the representation $R_N$ to a unitary representation
$\wt R_N$ of the group $U(N)\times U(N)$ acting in the same
space $L^2(U(N))$, cf. \tht{1.1}:
$$
(\wt R_N(g_1,g_2)f)(x)=f(g_2^{-1}xg_1), \qquad (g_1,g_2)\in
U(N)\times U(N).
$$
We call $\wt R_N$ the {\it biregular representation.}

In contrast to $R_N$, the decomposition of $\wt R_N$ is {\it
multiplicity free\/}:
$$
\wt R_N=\bigoplus_{\la\in\wh{U(N)}} (\pi^\la\otimes\pi^{\la^*}).
\tag2.1
$$
Here $\pi^{\la^*}$ stands for the conjugate representation to
$\pi^\la$; its signature is $\la^*=(-\la_N,\dots,-\la_1)$. We
observe that general irreducible representations of $U(N)\times
U(N)$ are of the form $\pi^\la\otimes\pi^\mu$, where
$\la,\mu\in\wh{U(N)}$. Representations with $\mu=\la^*$ are
characterized as those possessing a {\it spherical vector\/},
that is, a nonzero vector invariant under the subgroup
$\diag(U(N))$. Such representations are called {\it
spherical\/}. The whole subspace of $\diag(U(N))$--invariants in
$\pi^\la\otimes\pi^{\la^*}$ has dimension 1, so that the
spherical vector is defined uniquely up to a scalar factor.
Therefore, the spherical vector is a distinguished vector in the
representation space.

Note that the homogeneous space $U(N)\times U(N)/\diag(U(N))$ is
an example of a compact symmetric space $G/K$. For any such
space, the associated unitary representation of $G$ in
$L^2(G/K)$ is multiplicity free and its decomposition involves
exactly the irreducible spherical representations of the pair
$(G,K)$, that is, those irreducible representations of $G$ that
possess a $K$--invariant vector.

Returning to our special situation we conclude that the dual
space $\wh{U(N)}$ admits an alternative interpretation as the
set of (equivalence classes of) irreducible spherical
representations of the pair $(G,K)=(U(N)\times U(N),
\diag(U(N)))$.

Now we shall explain how this picture transforms when $U(N)$ is
replaced by $U(\infty)$.

\head 3. The dual space $\wh{U(\infty)}$ and spherical
representations of $U(\infty)\times U(\infty)$
\endhead

Consider the tower of groups $U(1)\subset U(2)\subset
U(3)\subset\dots$ where, for each $N$, the group $U(N)$ is
identified with the subgroup in $U(N+1)$ formed by matrices
$g=[g_{ij}]$ such that $g_{i,n+1}=g_{n+1,i}=\de_{i,n+1}$. We
define $U(\infty)$ as the union of all groups $U(N)$.
Equivalently, $U(\infty)$ consists of unitary matrices
$g=[g_{ij}]$ of infinite size, such that $g_{ij}=\de_{ij}$ for
$i+j$ large enough.

The conventional definition of a dual space, when applied to the
group $U(\infty)$, gives a huge pathological space.
\footnote{This is a general property of the so-called {\it wild
groups\/}; $U(\infty)$ is one of them.} It turns out that the
situation drastically changes if we mimic the alternative
interpretation of $\wh{U(N)}$ stated at the end of \S2:

\example{Definition 3.1} We set $\wh{U(\infty)}$ to be the space
of (equivalence classes of) irreducible spherical unitary
representations of the pair $(G,K)$, where
$$
G=U(\infty)\times U(\infty), \qquad K=\diag(U(\infty)). \tag3.1
$$
\endexample

Here ``spherical'' has the same meaning as above: existence of a
nonzero $K$--invariant vector. Again, such a vector is then
unique, within a scalar factor. Below $\R_+\subset\R$ denotes
the set of nonnegative real numbers and $\R_+^\infty$ denotes
the direct product of countably many copies of $\R_+$.

\proclaim{Theorem 3.2} The space $\wh{U(\infty)}$, see
Definition 3.1, can be identified with the subset
$\Om\subset\R_+^{4\infty+2}=\R_+^\infty\times\R_+^\infty
\times\R_+^\infty\times\R_+^\infty\times\R_+\times\R_+$ formed
by 6--tuples
$\om=(\al^+,\,\be^+,\,\al^-,\,\be^-,\,\de^+,\,\de^-)$ such that
$$
\gather
\al^{\pm}=(\al^\pm_1\ge\al^\pm_2\ge\dots\ge0)\in\R^\infty_+,
\quad
\be^{\pm}=(\be^\pm_1\ge\be^\pm_2\ge\dots\ge0)\in\R^\infty_+,\\
\de^\pm\in\R_+, \qquad \be^+_1+\be^-_1\le1, \quad
\sum_{i\ge1}(\al^\pm_i+\be^\pm_i)\le\de^\pm.
\endgather
$$
\endproclaim

Thus, for any point $\om\in\Om$ there exists an attached
irreducible spherical representation of $(G,K)$ which we denote
by $T^\om$. Representations $T^\om$ enter a larger class of {\it
admissible representations\/} which are studied in detail in
Olshanski \cite{Ol5}, \cite{Ol3}. In particular, we dispose of
an explicit description of the representation space of $T^\om$
together with the action of $G$ in it.

Theorem 3.2 has a long history. First of all, it should be said
that the classification of irreducible spherical representations
of $(G,K)$ is equivalent to that of {\it finite factor
representations\/} of the group $U(\infty)$, see Olshanski
\cite{Ol1}, \cite{Ol5, \S24}. \footnote{About factor
representations, see, e.g., Naimark \cite{Na, \S41.5}. In the
present paper we do not use this concept.} Finite factor
representations of $U(\infty)$ were first studied by Voiculescu
\cite{Vo}. He discovered (among many other things) that these
representations are parametrized by the so--called two--sided
infinite totally positive sequences of real numbers. But he did
not know that such sequences were completely classified much
earlier by Edrei \cite{Ed}. This fact was pointed out later by
Vershik--Kerov \cite{VK2} and Boyer \cite{Boy}. Thus, Theorem
3.2 is hidden in Edrei's paper. Note that \cite{Ed} is a pure
analytical work, which at first glance has nothing in common
with representation theory. Another, very different approach to
Theorem 3.2 was suggested in Vershik--Kerov \cite{VK2} and
further developed in Okounkov-Olshanski \cite{OkOl}.

Let $\SGN(N)\subset\Z^N$ denote the set of signatures of length
$N$, see \S2. We shall now define a sequence of embeddings
$\iota_N:\SGN(N)\to\Om$ such that as $N\to\infty$, the image
$\iota_N(\SGN(N))$ becomes more and more dense in $\Om$. This
agrees with the intuitive idea that the space $\wh{U(\infty)}$
should be a limit (in an appropriate sense) of the spaces
$\wh{U(N)}$. First, we need

\example{Definition 3.3 (Vershik--Kerov \cite{VK1})} Let $\mu$
be a Young diagram, $\mu'$ denote the transposed diagram, and
$d(\mu)$ denote the number of diagonal boxes in $\mu$. We also
regard $\mu$ as a partition $\mu=(\mu_1,\mu_2,\dots)$, so that
$\mu_i$ is the length of the $i$th row in $\mu$ while $\mu'_i$
is the length of the $i$th column. The numbers
$$
a_i(\mu)=\mu_i-i+\tfrac12, \quad b_i(\mu)=\mu'_i-i+\tfrac12,
\qquad 1\le i\le d(\mu)
$$
are called the {\it modified Frobenius coordinates\/} of $\mu$.
\endexample

For instance, if $\mu$ is the partition $(3,3,1,0,0,\dots)$ then
$d(\mu)=2$ and $a_1(\mu)=2\tfrac12$, $a_2(\mu)=1\tfrac12$,
$b_1(\mu)=2\tfrac12$, $b_2(\mu)=\tfrac12$. The modified
Frobenius coordinates are always positive half-integers whose
sum equals $|\mu|$, the number of boxes in $\mu$.

\example{Definition 3.4 (Embedding $\iota_N:\SGN(N)\to\Om$)}
Given a signature $\la\in\SGN(N)$, we represent it as a couple
$(\la^+,\la^-)$ of Young diagrams corresponding to positive and
negative coordinates in $\la$:
$$
\la=(\la^+_1\ge\dots\ge\la^+_k>0,\dots,0>-\la^-_l\ge\dots\ge-\la^-_1).
$$
Then we assign to $\la$ a point $\om=\iota_N(\la)\in\Om$, see
Theorem 3.2, as follows
$$
\al^\pm_i=\cases \frac{a_i(\la^\pm)}{N}, & i\le d(\la^\pm)\\ 0,
& i>d(\la^\pm) \endcases\,; \quad \be^\pm_i=\cases
\frac{b_i(\la^\pm)}{N}, & i\le d(\la^\pm)\\ 0, &
i>d(\la^\pm)\endcases\,; \quad \de^\pm=\frac{|\la^\pm|}{N}\,.
$$
\endexample

It is readily verified that $\om=(\al^+,\be^+,
\al^-,\be^-,\de^+,\de^-)$ is indeed a point of $\Om$. In
particular, the inequality $\be^+_1+\be^-_1\le1$ follows from
the evident fact that $k+l\le N$.

We equip $\Om$ with the topology inherited from the ambient
product space $\R_+^{4\infty+2}$. Then any point $\om\in\Om$ can
be approached by a sequence of the form $\iota_N(\la^{(N)})$,
where $\la^{(N)}\in\SGN(N)$, $N\to\infty$. Moreover, given a
sequence $\{\la^{(N)}\}$, we have
$$
\left(\iota_N(\la^{(N)})\to\om\right) \quad\Leftrightarrow\quad
\left(\pi^{\la^{(N)}}\otimes\pi^{{\la^{(N)}}^*}\to T^\om\right),
$$
where the last arrow means the convergence of representations of
the groups $U(N)\times U(N)$ to a representation of the group
$G=U(\infty)\times U(\infty)$, as defined in Olshanski
\cite{Ol5, \S22}, \cite{Ol2}.

\head 4. The problem of harmonic analysis \endhead

Let us try to understand now what could be an analog of the
decomposition \tht{2.1} for the group $G$. {}From \S3 we already
know the counterparts of the discrete set $\wh{U(N)}$ and of the
representations $\pi^{\la}\otimes\pi^{\la^*}$: these are the
infinite--dimensional space $\Om$ and spherical representations
$T^\om$. But what is the counterpart of the biregular
representation $\wt R_N$ acting in the Hilbert space
$L^2(U(N))$?

The conventional definition is not applicable to the group
$U(\infty)$: one cannot define the $L^2$ space on this group,
because $U(\infty)$ is not locally compact and hence does not
possess an invariant measure. To surpass this difficulty we
embed $U(\infty)$ into a larger space $\frak U$, which can be
defined as a {\it projective limit\/} of the spaces $U(N)$ as
$N\to\infty$. The space $\frak U$ is no longer a group but it is
still a $G$--space. That is, the two--sided action of
$U(\infty)$ on itself can be extended to an action on the space
$\frak U$. In contrast to $U(\infty)$, the space $\frak U$
possesses a biinvariant finite measure, which should be viewed
as a substitute of the nonexisting Haar measure. Moreover, this
biinvariant measure is included into a whole family
$\{\mu^{(s)}\}_{s\in\C}$ of measures with good transformation
properties. \footnote{The idea to enlarge an
infinite--dimensional space in order to build measures with good
transformation properties is well known. This is a standard
device in measure theory on linear spaces, but there are not so
many works where it is applied to ``curved'' spaces (see,
however, Pickrell \cite{Pi1}, Neretin \cite{Ner}). For the
history of the measures $\mu^{(s)}$ we refer to Olshanski
\cite{Ol7} and Borodin--Olshanski \cite{BO5}. A parallel
construction for the symmetric group case is given in
Kerov--Olshanski--Vershik \cite{KOV1}, \cite{KOV2}.} Using the
measures $\mu^{(s)}$ we explicitly construct a family
$\{T_{z,w}\}_{z,w\in\C}$ of representations, which seem to be a
good substitute of the nonexisting biregular representation. In
our understanding, the $T_{z,w}$'s are ``natural
representations'', and we state the problem of harmonic analysis
on $U(\infty)$ as follows:

\proclaim{Problem 4.1} Decompose the representations $T_{z,w}$
on irreducible representations.
\endproclaim

We skip a concrete description of the representations $T_{z,w}$,
which can be found in Olshanski \cite{Ol7}, and only list some
of their properties that are relevant for our discussion.
Henceforth we will assume that $\Re(z+w)>-1$ and that $z$ and
$w$ are not integers. Then, as it follows from the construction,
$T_{z,w}$ comes with a distinguished unit vector $\xi$, which is
$K$--invariant and cyclic. The latter property means that the
linear span of the $G$--orbit of $\xi$ is dense in
$H=H(T_{z,w})$, the Hilbert space of $T_{z,w}$. Let $H_N\subset
H$ be the Hilbert subspace spanned by the orbit of $\xi$ under
the subgroup $U(N)\times U(N)\subset G$. Then $H_N$ carries a
unitary representation of $U(N)\times U(N)$, which turns out to
be equivalent to the biregular representation $\wt R_N$ of \S2.
Since $\{H_N\}$ is an ascending chain of spaces whose union is
dense in $H$, we see that {\it $T_{z,w}$ is an inductive limit
of the biregular representations $\wt R_N$.\/} At this place the
reader might ask about the meaning of parameters $z,w$; the
answer is that to each value of $(z,w)$ there corresponds a
specific tower of embeddings
$$
H_1=L^2(U(1))\subset\dots\subset H_N=L^2(U(N)) \subset
H_{N+1}=L^2(U(N+1))\subset\dots. \tag4.1
$$
There are many (even too many) ways to realize $\wt R_N$ as a
subrepresentation of $\wt R_{N+1}$, and our construction leads
to a distinguished 2--parameter family of towers of embeddings.

The statement of Problem 4.1 looks rather abstract but we will
gradually reduce it to a concrete form. The first step is to
apply the following abstract claim.

\proclaim{Theorem 4.2} Let $T$ be a unitary representation of
$G$ in a Hilbert space $H$ and assume that there exists  a
$K$--invariant cyclic vector $\xi\in H$ (we will assume
$\Vert\xi\Vert=1$). Then $(T,\xi)$ is completely determined,
within a natural equivalence, by a probability measure $P$ on
the dual space $\wh{U(\infty)}=\Om$. The decomposition of $T$ on
irreducible representations is given by a multiplicity free
direct integral of spherical representations $T^\om$ with
respect to measure $P$.
\endproclaim

We call $P$ the {\it spectral measure\/} of $(T,\xi)$. Note that
if $\xi$ is replaced by another vector $\xi'\in H$ with the same
properties then $P$ is replaced by an equivalent measure $P'$.
We will not define precisely what is a ``direct integral of
representations'' (see, e.g., Naimark \cite{Na, \S41}) but only
observe that Theorem 4.2 is strictly similar to a customary
fact, the spectral theorem for a pair $(A,\xi)$ where $A$ stands
for a self--adjoint operator in a Hilbert space $H$ and $\xi\in
H$ is a unit cyclic vector.

Taking into account Theorem 4.2 we replace Problem 4.1 by

\proclaim{Problem 4.3} Assume that $z,w\in\C\setminus\Z$ and
$\Re(z+w)>-1$. Let $\xi$ be the distinguished $K$--invariant
cyclic unit vector provided by the construction of $T_{z,w}$,
and let $P_{z,w}$ denote the spectral measure of
$(T_{z,w},\xi)$, which is a probability measure on $\Om$.
Describe $P_{z,w}$ explicitly.
\endproclaim

Recall that the Hilbert space $H(T_{z,w})$ is the inductive
limit of a chain \tht{4.1} and that the vector $\xi$ belongs to
all spaces $H_N$, which carry representations $\wt R_N$.
Evidently, for each $N$, $\xi$ is a $\diag(U(N))$--invariant
cyclic vector in the biregular representation $\wt R_N$. The
pair $(\wt R_N,\xi)$ gives rise to a spectral measure
$P^{(N)}_{z,w}$ on $\wh{U(N)}=\SGN(N)$. Since $\SGN(N)$ is a
discrete space, this is a purely atomic probability measure. It
has a very simple meaning. According to decomposition \tht{2.1}
we obtain an orthogonal decomposition of $\xi$ into a sum of
certain vectors $\xi_\la$. We have
$$
1=\Vert\xi\Vert^2=\sum_{\la\in\SGN(N)}\Vert\xi_\la\Vert^2 \quad
\text{and}\quad P^{(N)}_{z,w}(\la)=\Vert\xi_\la\Vert^2 \quad
\text{for} \quad \la\in\SGN(N).
$$

The numbers $P^{(N)}_{z,w}(\la)$ can be computed, the result is
as follows
$$
\gather P^{(N)}_{z,w}(\la)=\const_N\cdot\prod_{i=1}^N
W_N(\la_i-i)\cdot
\prod_{1\le i<j\le N}(\la_i-\la_j-i+j)^2, \tag4.2\\
W_N(l)=\left|\Ga(z-l)\Ga(w+N+1+l)\right|^{-2}, \quad l\in\Z,
\tag4.3 \endgather
$$
where $\const_N$ is a normalization constant. The assumption
that $z,w$ are not integers just means that $P^{(N)}_{z,w}(\la)$
does not vanish (which is related to cyclicity of vector $\xi$).
The assumption $\Re(z+w)>-1$ guarantees that
$$
\sum_{\la\in\SGN(N)}\; \prod_{i=1}^N W_N(\la_i-i)\, \prod_{1\le
i<j\le N}(\la_i-\la_j-i+j)^2<\infty
$$
for all $N$, so that the normalization is indeed possible.

On the other hand, one can prove that
$$
\lim_{N\to\infty}\iota_N\left(P^{(N)}_{z,w}\right)=P_{z,w}\,,
\tag4.4
$$
where the embeddings $\iota_N:\SGN(N)\to\Om$ were specified in
Definition 3.4. Thus,  Problem 4.3 admits a reformulation which
already has a very concrete form:

\proclaim{Problem 4.4} Compute explicitly the limit probability
measure in the right--hand side of \tht{4.3}, where the
probability measures in the left--hand side are given by
\tht{4.2} and Definition 3.4.
\endproclaim

In the remaining part of the paper we explain how this problem
is solved. A detailed exposition of the material of this section
can be found in Olshanski \cite{Ol7}.

\head 5. Random point processes \endhead

The spectral measures $P_{z,w}$ that we aim to describe live on
a ``very big'' space $\Om$, which is a domain in an
infinite--dimensional product space. There is no hope that $\Om$
possesses a simple reference measure (like Lebesgue measure)
such that $P_{z,w}$ would be determined by a density with
respect to that measure. Thus, we have to use another language
to describe our measures. It turns out that such a language is
provided by the theory of random point processes.

In this section we give a few necessary basic definitions
concerning random point processes and also provide a few
examples which seem to be relevant for the discussion of our
main problem. One should not regard our exposition as a survey
on point processes. As basic references on this subject the
reader can consult Daley and Vere-Jones \cite{DVJ} and Lenard
\cite{Len}

Let $\X$ be a locally compact space. A {\it point
configuration\/} in $\X$ is a finite or countable subset without
limit points. Let $\Conf(\X)$ be the set of all point
configurations. For any Borel subset $A\subset\X$ with compact
closure, let $\N_A:\Conf(\X)\to\Z_+$ be the function defined by
$\N_A(X)=|A\cap X|$, where $X\in\Conf(\X)$. Consider the
sigma--algebra of subsets in $\Conf(\X)$ generated by all
functions $\N_A$. A probability measure $\Cal P$ defined on this
sigma--algebra is called a {\it random point process\/} on $\X$.
Given $\Cal P$, point configurations $X\subset\X$ become random
objects, and we can speak, for instance, about probabilities of
events like this:
$$
\N_{A_1}(X)=n_1,\dots,\N_{A_k}(X)=n_k.
$$

\example{Example 5.1 (Poisson process)} The simplest and most
known random point process is the Poisson process, which is
determined by an arbitrary measure $m$ on $\X$. The Poisson
process is characterized by the property that the probability of
each event of the form above, where $A_1,\dots,A_k$ do not
intersect, equals
$$
\prod_{i=1}^k e^{-m(A_i)}\frac{(m(A_i))^{n_i}}{n_i!}\,.
$$
\endexample

Given a point process $\Cal P$ on $\X$, we can integrate various
functions $F(X)$ on $\Conf(\X)$. An important class of functions
$F$ is defined as follows. Let $f(x_1,\dots,x_n)$ be a
continuous function on $\X^n$ with compact support; we set
$$
F_f(X)=\sum_{x_1,\dots,x_n}f(x_1,\dots,x_n), \qquad
X\in\Conf(\X),
$$
summed over all $n$--tuples of pairwise distinct points in $X$.
Note that $F_f$ depends on the symmetric part of $f$ only. Under
mild assumptions on $\Cal P$, there exists a unique symmetric
measure $\rho_n$ on $\X^n$ such that for any $f$ as above,
$$
\int\limits_{\Conf(\X)}F_f(X)\Cal P(dX)
=\int\limits_{\X^n}f(x_1,\dots x_n)\rho_n(dx_1\dots dx_n),
$$
and, moreover, $\Cal P$ is uniquely determined by the infinite
sequence of measures $\rho_1,\rho_2,\dots$ (see Lenard
\cite{Len}). These measures are called the {\it correlation
measures\/} of $\Cal P$. They are a convenient tool for
identifying and studying a point process.

When $\Cal P$ is the Poisson process, we simply have
$\rho_n=m^{\otimes n}$. For non--Poisson processes $\Cal P$, the
correlation measures can have a more sophisticated structure.

In practice one can usually choose a natural reference measure
$m$ on $\X$ such that $\rho_n$ has a density with respect to
$m^{\otimes n}$  for each $n$. Then this density is called the
$n$th {\it correlation function\/} of $\Cal P$; we will denote
it as $\rho_n(x_1,\dots,x_n)$. If $\X$ is a discrete space and
$m$ is the counting measure then $\rho_n(x_1,\dots,x_n)$ is the
probability that the random configuration $X$ contains all
points $x_1,\dots,x_n$ (if these points are not all distinct
then $\rho_n(x_1,\dots,x_n)=0$). When $\X$ is not discrete,
$\rho_n(x_1,\dots,x_n)$ can be informally defined as follows
$$
\rho_n(x_1,\dots,x_n)=\lim_{\De x_1\to0, \dots,\De x_n\to0}
\frac{\text{$\Prob\{$ random $X$ intersects $\De
x_1$,\dots,\,$\De x_n\}$}}{m(\De x_1)\dots m(\De x_n)}\,,
$$
where $\De x_1,\dots,\De x_n$ are small neighborhoods of the
points $x_1,\dots,x_n$. In words, $\rho_n(x_1,\dots,x_n)$ is the
density of the probability to find a point of the random
configuration in each of $n$ infinitesimally small neighborhoods
about $x_1,\dots,x_n$.

\example{Definition 5.2 (Determinantal processes)} Assume that a
reference measure as above exists, so that we can deal with the
correlation functions. Then $\Cal P$ is called a {\it
determinantal point process\/} if there exists a function
$K(x,y)$ on $\X\times\X$ such that
$$
\rho_n(x_1,\dots,x_n)=\det[K(x_i,x_j)]_{1\le i,j\le n}\,, \qquad
n=1,2,\dots.
$$
We call $K$ the {\it correlation kernel\/} of $\Cal P$.
\endexample

If $K$ is symmetric ($K(x,y)=\overline{K(y,x)}$) then the points
in the random configuration are {\it negatively correlated\/}: a
very close rapprochement of points has a relatively small
probability. So, the points look as mutually repelling
particles. In a Poisson process, on the contrary, the points are
not correlated at all; they look as noninteracting particles. A
good survey on determinantal point processes is Soshnikov
\cite{So}.

All the information about a determinantal process $\Cal P$ is
hidden in its correlation kernel $K(x,y)$. In this respect,
determinantal point processes can be compared to Gaussian
measures where all the information is contained in the
covariation matrix. Knowing $K(x,y)$ we can, in principle,
compute the probabilities of various natural events associated
to $\Cal P$. We state the simplest but important example:

\proclaim{Proposition 5.3} Let $\Cal P$ be a determinantal point
process with a correlation kernel $K$. The probability of having
no particles in a region $I\subset\X$ is equal to the Fredholm
determinant $\det(1-K_I)$, where $K_I$ is the restriction of $K$
to $I\times I$.
\endproclaim

It often happens that such {\it gap probability} can be
expressed through a solution of a (second order nonlinear
ordinary differential) Painlev\'e equation, see Example 6.2
below.

The most known example of a determinantal process is

\example{Example 5.4 (Sine process)} The {\it sine kernel\/} is
given by
$$
K(x,y)=\frac{\sin(\pi(x-y))}{\pi(x-y)}\,, \qquad x,y\in\R
$$
(here the reference measure $m$ is Lebesgue measure). The sine
kernel determines a remarkable translation invariant point
process on $\X=\R$.
\endexample

It is instructive to compare the sine process with the standard
Poisson process on $\R$ (where $m$ is again Lebesgue measure).
Both processes are translation invariant, and for both processes
the mean distance between adjacent points equals 1. However, as
can be seen from computer simulations, the sample random
configurations of the Poisson process are more chaotic. For the
Poisson process, the distance between adjacent points is a very
simple random variable (it has exponential distribution), while
for the sine process the corresponding distribution is expressed
through a Painlev\'e transcendent.\footnote{This result was
originally proved in Jimbo--Miwa-- M\^ori--Sato\cite{JMMS}, and
a number of other proofs and extensions were later given by
different authors, see Borodin--Deift \cite{BD} for
references.}.

For a large number of concrete examples of determinantal
processes the space $\X$ is a subset of $\R$, $\C$, or $\Z$, and
the correlation kernel has the form
$$
K(x,y)=\frac{P(x)Q(y)-Q(x)P(y)}{x-y} \tag5.1
$$
or, more generally,
$$
K(x,y)=\frac{\sum\limits_{i=1}^k F_i(x)G_i(y)}{x-y}\,, \qquad
\text{where}\quad \sum\limits_{i=1}^kF_i(x)G_i(x)=0.\tag5.2
$$
Such kernel are called {\it integrable\/}, see
Its--Izergin--Korepin--Slavnov \cite{IIKS}, Deift \cite{De},
Borodin \cite{B2}.

\example{Example 5.4 (Orthogonal polynomial ensembles)} Let
$W(x)$ be a weight function (defined, say, on a subset
$\X\subset\R$) and let $p_0\equiv1, p_1, p_2,\dots$ be the
associated family of orthogonal polynomials. For an arbitrary
$N=1,2\dots$, consider the orthogonal projection operator in
$L^2(\X,dx)$ \,\footnote{If $\X$ is a discrete set then Lebesgue
measure $dx$ is replaced by the counting measure.} onto the
$N$--dimensional subspace spanned by functions
$p_i(x)W^{\frac12}(x)$, $0\le i\le N-1$, and let $K^W_N(x,y)$
stand for the kernel of this operator. This kernel can be
written in integrable form \tht{5.1} with
$$
P(x)=\const\, p_N(x)W^{\frac12}(x), \qquad Q(x)=\const\,
p_{N-1}(x)W^{\frac12}(x).
$$
In other words, $K^W_N(x,y)$ is equal to the classical
Christoffel--Darboux kernel times
$W^{\frac12}(x)W^{\frac12}(y)$. The kernel $K^W_N(x,y)$ gives
rise to random $N$--point configurations in $\X$. Namely, the
density of probability \,\footnote{If the space $\X$ is discrete
then one can simply speak about the probability of
$(x_1,\dots,x_N)$.} of a given configuration has the form
$$
\Cal P(x_1,\dots,x_N)=\const\,\prod_{i=1}^N W(x_i)\,\prod_{1\le
i<j\le N} (x_i-x_j)^2. \tag5.3
$$
The random point processes of this type are called {\it
orthogonal polynomial ensembles\/.} Note that \tht{5.3} can be
written in the Gibbsian form which is common in statistical
physics:
$$
\Cal P(x_1,\dots,x_N)=\const\,\exp\left(-\sum_i\log
V^{-1}(x_i)-2\sum_{i<j}\log|x_i-x_j|^{-1}\right).
$$
The terms $\log V^{-1}(x)$ and $2\log|x_i-x_j|^{-1}$ are
interpreted as the one--particle potential and the pair
potential, respectively, and the whole ensemble is interpreted
as an $N$--particle {\it log--gas system\/} (Forrester
\cite{Fo}).

\endexample

A variety of random point processes comes from spectra of random
matrices. A basic example is the Gaussian Unitary Ensemble (GUE)
formed by $N\times N$ Hermitian matrices distributed according
to a Gaussian measure invariant under conjugation by unitary
matrices from $U(N)$. The spectrum of such a random matrix is a
random $N$--point configuration in $\X=\R$ arising from the
Hermite orthogonal polynomial ensemble (in the notation of
Example 5.4, $W(x)=e^{-x^2}$, the weight function of the Hermite
polynomials). {}From other ensembles of random matrices one can
also obtain the Laguerre and Jacobi  orthogonal polynomial
ensembles (see, e.g., Forrester \cite{Fo}).

One of the fundamental problems in random matrix theory is to
study the asymptotic behavior of random matrices as their size
goes to infinity. This leads, in particular, to studying the
{\it scaling limits\/} of orthogonal polynomial ensembles in
various regimes. For instance, if we focus at the $N$--point
Hermite polynomial ensemble with large $N$ in a neighborhood of
the origin and scale the space variable $x$ so that the mean
distance between adjacent points becomes approximately 1 (which
is achieved by the change of variable $x\to x'=\sqrt{2N}x/\pi$),
then we obtain in the limit $N\to\infty$ the sine process.

Orthogonal polynomial ensembles with discrete state space $\X$
arise in a number of probabilistic models which include random
tilings (Johansson \cite{Jo3}) and directed percolation
(Johansson \cite{Jo1}, \cite{Jo2}). Classical discrete
orthogonal polynomials known as Charlier, Krawtchouk, Meixner,
and Hahn polynomials arise in this fashion.

\head 6. Point processes $\Cal P_{z,w}$. The main result
\endhead

Now we return to the spectral measures $P_{z,w}$. We will
explain how to convert them into random point processes $\Cal
P_{z,w}$ on the space
$$
\X=\R\setminus\{\pm\tfrac12\}
$$
(the real line with two punctures, at $\frac12$ and $-\frac12$).

We define a projection $\Om\to\Conf(\X)$ by
$$
\multline
\om=(\al^+,\be^+,\al^-,\be^-,\de^+,\de^-) \\
\mapsto X=\{\al_i^++\tfrac{1}2\}\sqcup\{\tfrac{1}2-
\be_i^+\}\sqcup\{-\al_j^--\tfrac{1}2\}\sqcup\{-
\tfrac{1}2+\be_j^-\},
\endmultline \tag6.1
$$
where we omit possible 0's among
$\al^+_i,\be^+_i,\al^-_i,\be^-_i$, and also omit possible 1's
among $\be^+_i$ or $\be^-_i$. Note that $X$ is bounded in $\R$
and its points may accumulate only near the punctures $\frac12$
and $-\frac12$.

By definition, $\Cal P_{z,w}$ is the push--forward of the
measure $P_{z,w}$ under the projection $\Om\to\Conf(\X)$.

The projection is not injective, so that we can, in principle,
loose a part of information about our measure $P_{z,w}$ under
the passage $P_{z,w}\to\Cal P_{z,w}$. However, one can present
arguments showing that the losses (if any) are negligible, see
the end of \S9 in Borodin--Olshanski \cite{BO6}. Thus, we can
regard $\Cal P_{z,w}$ as a substitute of $P_{z,w}$.

The next result provides a description of the point process
$\Cal P_{z,w}$ and can be viewed as a solution of Problem 4.4.

\proclaim{Theorem 6.1 (Main result)} $\Cal P_{z,w}$ is a
determinantal point processes. Its correlation kernel can be
written in integrable form \tht{5.2} with $k=2$, where the
functions $F_1$, $F_2$, $G_1$, and $G_2$ can be explicitly
expressed through the Gauss hypergeometric function.

For instance, if $x>\frac 12$ and $y>\frac 12$ then the kernel
can be written in form \tht{5.1} with
$$
\multline P(x)=\const\,\left(x-\frac 12\right)^{-\frac12(z+\bar
z)-\bar w}\left(x+\frac
12\right)^{\frac12(\bar w- w)}\\
\times{}_2F_1\left(z+\bar w,\,\bar z+\bar w; \, z+\bar z+w+\bar
w+1;\,\left(\tfrac 12 -x\right)^{-1}\right),
\endmultline
$$
$$
\multline Q(x)=\const\,\left(x-\frac 12\right)^{-\frac12(z+\bar
z)-\bar w-1}\left(x+\frac
12\right)^{\frac12(\bar w-w)}\\
\times{}_2F_1\left(z+\bar w+1,\,\bar z+\bar w+1; \, z+\bar
z+w+\bar w+2;\,\left(\tfrac 12 -x\right)^{-1}\right).
\endmultline
$$
\endproclaim

Here ${}_2F_1(a,b;c;\zeta)$ is the Gauss hypergeometric function
with parameters $a,b,c$ and argument $\zeta$. Note that this
function is well defined for $\zeta<0$.

We call the kernel of Theorem 6.1 the (continuous) {\it
hypergeometric kernel\/}; let us denote it by
$K^\hypgeom_{z,w}(x,y)$. Precise formulas for the kernel and the
proof of the theorem are given in our paper \cite{BO6}.

Note that the kernel $K^\hypgeom_{z,w}(x,y)$ is real valued but
{\it not\/} symmetric. It has the following symmetry property
instead:
$$
K^\hypgeom_{z,w}(x,y)=\cases K^\hypgeom_{z,w}(y,x) &\text{if
$x,y$ are both
inside}\\ &\text{or outside $(-\frac12,\frac12)$;}\\
-K^\hypgeom_{z,w}(y,x) &\text{otherwise.}\endcases \tag6.2
$$
In other words, $K^\hypgeom_{z,w}(x,y)$ is symmetric with
respect to the indefinite inner product of functions on $\X$
given by
$$
[f,g]=\int\limits_{\R\setminus[-\frac12,\frac12]}f(x)g(x)dx
\quad -\quad \int\limits_{(-\frac12,\frac12)}f(x)g(x)dx
$$
An explanation of this fact will be given in Remark 7.2 below.

Since all the information about the point process $\Cal P_{z,w}$
is hidden in the kernel $K^\hypgeom_{z,w}(x,y)$, a natural
question is: What can be extracted from the explicit expression
for the kernel? For instance, each of parameters $\al^\pm_i$,
$\be^\pm_i$ can be viewed as a random variable defined on the
probability space $(\Om,P_{z,w})$; what can be said about their
distribution? Here are two examples.

The first example concerns the distribution of $\al^+_1$. The
same result holds for $\al^-_1$; it suffices to interchange $z$
and $w$.

\example{Example 6.2 (Painlev\'e VI)} By virtue of Proposition
5.3, the probability distribution of $\al^+_1$ is given by
$$
\Prob\{\al^+_1<u\}=\det(1-K_{\frac12+u}), \qquad u>0,
$$
where we abbreviate
$$
K_s=\left.K^\hypgeom\right|_{(s,+\infty)\times(s,+\infty)}\,,
\qquad s>\tfrac12\,.
$$
Set
$$
\gathered \nu_1=\frac{z+\bar z+w+\bar w}2\,,\quad
\nu_3=\frac{z-\bar z+w-\bar w}2\,,\quad
\nu_4=\frac{z-\bar z-w+\bar w}2\,,\\
\sigma(s)=\left(s^2-\tfrac 14\right)\frac {d\ln\det(1-K_s)}{ds}-
\nu_1^2\,s+\frac {\nu_3\nu_4}{2}.
\endgathered
$$
Then $\sigma(s)$ satisfies the differential equation
$$
-\sigma'\left(\left(s^2-\tfrac
14\right)\sigma''\right)^2=\left(2\left(s\sigma'-
\sigma\right)\sigma' -\nu_1^2\nu_3\nu_4\right)^2 -
(\sigma'+\nu_1^2)^2(\sigma'+\nu_3^2)(\sigma'+\nu_4^2).
$$
This differential equation is the so--called $\sigma$-form of
the Painlev\'e VI equation. The proof can be found in
Borodin--Deift \cite{BD}. We refer to the introduction of that
paper for a brief historical introduction and references on this
subject.
\endexample

Our second example concerns the asymptotic behavior of
parameters $\al^\pm_i$, $\be^\pm_i$ as $i\to\infty$.

\example{Example 6.3 (Law of large numbers)} We conjecture that
with probability 1,
$$
\lim_{k\to\infty}(\al^+_k)^{1/k}=\lim_{k\to\infty}(\be^+_k)^{1/k}=q(z),
\qquad
\lim_{k\to\infty}(\al^-_k)^{1/k}=\lim_{k\to\infty}(\be^-_k)^{1/k}=q(w),
$$
where
$$
q(z)=\exp\left(-\,\sum_{n\in\Z}|z-n|^{-2}\right)
=\exp\left(-\,\frac{\pi\sin(\pi(z-\bar z))}{(z-\bar z)\sin(\pi
z)\sin (\pi\bar z)}\right)
$$
This conjecture is based on the results of Borodin--Olshanski
\cite{BO1} and \cite{BO7}. The result should be obtained by
analogy with Theorem 5.1 of \cite{BO1}. However, we did not
verify the details yet.
\endexample

\head 7. Lattice approximation to process $\Cal P_{z,w}$
\endhead

Our proof of Theorem 6.1 is based on the limit relation
\tht{4.4}. In \S6, we have interpreted its right--hand side as a
point process. Here we explain how to do the same for the
left--hand side and thus to translate this relation into the
language of random point processes.

Comparing \tht{4.2}--\tht{4.3} with \tht{5.3} we see that the
measure $P^{(N)}_{z,w}$ on $\SGN(N)$ gives rise to a discrete
orthogonal polynomial ensemble on $\Z$ with weight function
\tht{4.3}. Here we have used the bijective correspondence
between diagrams $\la\in\SGN(N)$ and $N$--point configurations
$(l_1>\dots>l_N)$ on $\Z$ determined by relation $l_i=\la_i-i$.

Since the weight $W_N(l)$ from \tht{4.3} has a slow (polynomial)
decay at infinity,
$$
W_N(l)\sim|l|^{-2\Re(z+w)-2N}\,, \qquad l\to\pm\infty,
$$
it admits only finitely many orthogonal polynomials. However,
due to the assumption $\Re(z+w)>-1$, we have enough polynomials
to define the orthogonal polynomial ensemble for any $N$. We
call it the {\it Askey--Lesky\/} ensemble, because the
orthogonal polynomials in question were computed in Askey
\cite{As} and Lesky \cite{Les1}, \cite{Les2}. The Askey--Lesky
polynomials are relatives of the classical Hahn polynomials;
they are expressed through the value of the hypergeometric
series $_3F_2$ at 1. {}From the explicit expression of these
polynomials we obtain the corresponding correlation kernel
$K^\askey_N(x,y)$. The Askey--Lesky ensemble is an interesting
example of a discrete log--gas system (the particles are
confined to a lattice).

However, the Askey--Lesky ensemble is only an intermediate
object, we need to transform it further in order to visualize
the modified Frobenius coordinates of Young diagrams $\la^\pm$
(see Definitions 3.3 and 3.4).

The first step is rather simple, we shift the configuration
$(l_1,\dots,l_N)$  by $\frac{N+1}2$, so that the resulting
correspondence between signatures and $N$--point configurations
takes a more symmetric form
$$
\la\,\leftrightarrow\,\Cal L=\{\la_1+\tfrac{N-1}2,\;
\la_2+\tfrac{N-3}2,\; \dots,\;
\la_{N-1}-\tfrac{N-3}2,\;\la_N-\tfrac{N-1}2\}. \tag7.1
$$
The configuration $\Cal L$ lives on the lattice
$$
\X^{(N)}=\Z+\tfrac{N+1}2=\cases \Z, & \text{if $N$ is odd;}\\
\Z+\frac12 &\text{if $N$ is even.} \endcases
$$

The next step is less obvious. Let us divide the lattice
$\X^{(N)}$ into two parts, which will be denoted by
$\X^{(N)}_{\inr}$ and $\X^{(N)}_{\out}$:
$$
\gathered
\X^{(N)}_{\inr}=\left\{-\tfrac{N-1}2,-\tfrac{N-3}2,\dots,\tfrac{N-3}2,
\tfrac{N-1}2\right\},\\
\X^{(N)}_{\out}=\left\{\dots,-\tfrac{N+3}2,-\tfrac{N+1}2\right\}
\cup\left\{\tfrac{N+1}2,\tfrac{N+3}2,\dots\right\}.
\endgathered
$$
Here $\X^{(N)}_{\inr}$, the ``inner'' part, consists of $N$
points of the lattice that lie on the interval $(-\tfrac
N2,\tfrac N2)$, while $\X^{(N)}_{\out}$, the ``outer'' part, is
its complement in $\X^{(N)}$, consisting of the points outside
this interval.

Given an $N$--point configuration $\Cal L$ on $\X^{(N)}$, which
we interpret as a system of particles occupying $N$ positions on
the lattice $\X^{(N)}$, we assign to it another configuration,
$X$, formed by the particles in $\X^{(N)}_{\out}$ and the {\it
holes\/} (i.e., the unoccupied positions) in $\X^{(N)}_{\inr}$.
Note that $X$ is a finite configuration, too. Since the
``interior'' part consists of exactly $N$ points, we see that in
$X$, there are equally many particles and holes. However, their
number is no longer fixed, it varies between 0 and $2N$,
depending on the mutual location of $\Cal L$ and
$\X^{(N)}_{\inr}$. For instance, if these two sets coincide then
$X$ is the empty configuration, and if they do not intersect
then $|X|=2N$.

We call the procedure of passage $\Cal L\mapsto X$ the {\it
particles/holes involution\/}. Under this procedure, our initial
random $N$--particle system (coming from the Askey--Lesky
ensemble) turns into a random system of particles and holes.
Note that the map $\Cal L\mapsto X$ is reversible, so that both
random point processes are equivalent. Let us denote the second
point process by $\Cal P^{(N)}_{z,w}$.

The significance of the procedure described above becomes clear
from the following combinatorial fact.

\proclaim{Lemma 7.1 (\cite{BO6, \S4})} Let $\la\in\SGN(N)$ be a
signature, $\Cal L\subset\X^{(N)}$ be the $N$--particle
configuration defined by \tht{7.1}, and $X\subset\X^{(N)}$ be
the corresponding finite configuration of particles and holes as
defined above. Let also $a^\pm_i$ and $b^\pm_i$ be the modified
Frobenius coordinates of the Young diagrams $\la^\pm$, see
Definitions 3.3 and 3.4.

Then we have
$$
X\cap\X^{(N)}_{\out}=\{a^+_i+\tfrac N2\}\cup\{-a^-_i-\tfrac
N2\}, \qquad X\cap\X^{(N)}_{\inr}=\{\tfrac
N2-b^+_i\}\cup\{-\tfrac N2+b^-_i\}.\tag7.2
$$
\endproclaim

Comparing \tht{7.2} with \tht{6.1} suggests that if we shrink
our phase space $\X^{(N)}$ by the factor of $N$ (so that the
points $\pm \frac N2$ turn into $\pm\frac 12$) then our discrete
point process $\Cal P^{(N)}_{z,w}$ should have a well--defined
scaling limit. We prove that such a limit does exist and it
coincides with the point process $\Cal P_{z,w}$  on $\X=\Bbb
R\setminus\{\pm \frac 12\}$ as defined in \S6.

The discrete process $\Cal P^{(N)}_{z,w}$ is determinantal, and
its correlation kernel can be obtained by a transformation of
the kernel $K^\askey_N(x,y)$; let us denote this new kernel by
$\wt K^\askey_N(x,y)$. The correlation kernel $K^\hypgeom(x,y)$
of Theorem 6.1 is obtained as a scaling limit of the kernel $\wt
K^\askey_N(x,y)$.

We just gave a rough sketch of the proof of Theorem 6.1. The
detailed proof (see Borodin--Olshanski \cite{BO6}) is rather
long and technical. The main technical difficulties arise when
we want to get a convenient explicit expression for the kernel
$\wt K^\askey_N(x,y)$ in case when at least one of variables
$x,y$ is in the ``interior'' part of the lattice. \footnote{This
part of the kernel describes the correlations of holes with
particles and other holes. The correlations involving particles
only are described by the kernel $K^\askey_N(x,y)$ restricted to
the ``exterior'' part of the lattice.} Here we apply a discrete
version of the formalism of the Riemann--Hilbert problem, see
Borodin \cite{B2}.

\example{Remark 7.2 (On symmetry \tht{6.2})} Now we are in a
position to explain the indefinite--type symmetry \tht{6.2}: the
same kind of symmetry occurs already in the kernel $\wt
K^\askey_N(x,y)$. It turns out that the particles/holes
involution just converts the usual symmetry of kernel
$K^\askey_N(x,y)$ into the indefinite--type symmetry of kernel
$\wt K^\askey_N(x,y)$.
\endexample

The point process $\Cal P^{(N)}_{z,w}$ can be viewed as a {\it
discrete two--component log--gas system\/} consisting of
oppositely signed charges. Systems of such a type were earlier
investigated in the mathematical physics literature (see, e.g.,
a number of references listed in section (f) of the introduction
to Borodin--Olshanski \cite{BO6}). However, the known concrete
models are quite different from our system.

\example{Remark 7.3 (Limit density)} Given an $N$--point
orthogonal polynomial ensemble, let us attach to a configuration
$\{x_1,\dots,x_N\}$ a probability measure,
$$
\tfrac1N(\de_{x_1}+\dots+\de_{x_N}).
$$
Under an appropriate scaling limit as $N\to\infty$, this random
measure can converge to a (nonrandom) probability measure
describing the global limit density of particles. For instance,
in case of GUE, the limit density is given by the famous
Wigner's semi--circle law, see e.g., Forrester \cite{Fo, ch. 1}.

When we apply this procedure to the Askey--Lesky ensemble (or
rather to its shift by $\frac{N+1}2$) then it can be shown that,
as $N$ gets large, almost all $N$ particles occupy positions
inside $(-\frac N2, \frac N2)$. (Recall that there are exactly
$N$ lattice points in this interval, hence, almost all of them
are occupied by particles.) In other words, this means that the
density of our discrete log--gas is asymptotically equal to the
characteristic function of the $N$--point set of lattice points
inside $(-\frac N2,\frac N2)$, so that in the scaling limit we
get the characteristic function of $(-\frac12,\frac12)$.

It can also be shown that after the passage $\Cal L\to X$, all
but finitely many particles/holes in $X$ concentrate, for large
$N$, near the points $\pm \frac N2$. This explains why the
random system of paricles/holes $X$ converges to a limit point
process (as opposed to the Askey--Lesky ensemble).
\endexample

\head 8. Connection with previous work \endhead

Let us briefly discuss  two similar problems  which also lead to
spectral measures on infinite--dimensional spaces.

The first problem was initially formulated in
Kerov--Olshanski--Vershik \cite{KOV1}. It consists in
decomposing certain natural (generalized regular) unitary
representations $T_z$ of the group $S(\infty)\times S(\infty)$,
depending on a complex parameter $z$. In \cite{KOV1},
\cite{KOV2} the problem was solved in the case when the
parameter $z$ takes integral values (then the spectral measure
has a finite--dimensional support). The general case presents
more difficulties and we studied it in a cycle of papers (see
the surveys Borodin--Olshanski \cite{BO2}, Olshanski \cite{Ol6}
and references therein). Our main result is that the spectral
measure governing the decomposition of $T_z$ can be described in
terms of a determinantal point process on the real line with one
punctured point. The correlation kernel was explicitly computed,
it has integrable form \tht{5.2}, where $k=2$ and the functions
$F_1$, $F_2$, $G_1$, and $G_2$ are expressed through a confluent
hypergeometric function (specifically, through the W--Whittaker
function), see Borodin \cite{B1}, Borodin--Olshanski \cite{BO3}.

The second problem deals with decomposition of a family of
unitarily invariant probability measures on the space of all
infinite Hermitian matrices on ergodic components. The measures
depend on one complex parameter; within a transformation of the
underlying space, they coincide with the measures $\mu^{(s)}$
mentioned in the beginning of \S4. The problem of decomposition
on ergodic components can be also viewed as a problem of
harmonic analysis on an infinite--dimensional Cartan motion
group. The main result states that the spectral measures in this
case can be interpreted as determinantal point processes on the
real line with an integrable correlation kernel of type
\tht{5.1}, where the functions $P$ and $Q$ are expressed through
another confluent hypergeometric function, the M--Whittaker
function, see Borodin--Olshanski \cite{BO5}.

These two problems and the problem that we deal with in this
paper have many similarities but the latter problem is, in a
certain sense, more general comparing to both problems described
above. The Askey--Lesky kernel of \S7 can be viewed as the top
of a hierarchy of (discrete and continuous) integrable kernels:
this looks very much like the hierarchy of the classical special
functions. A description of the ``$S(\infty)$--part'' of the
hierarchy can be found in Borodin--Olshanski \cite{BO4}.

\bigskip

{\smc A.~Borodin}: Mathematics 253-37, Caltech, Pasadena, CA
91125, U.S.A.,

\medskip

E-mail address: {\tt borodin\@caltech.edu}

\bigskip

{\smc G.~Olshanski}: Dobrushin Mathematics Laboratory, Institute
for Information Transmission Problems, Bolshoy Karetny 19,
127994 Moscow GSP-4, RUSSIA.

\medskip

E-mail address: {\tt olsh\@online.ru}

\enddocument

\Refs

\widestnumber\key{JMMS}

\ref\key As \by R.~Askey \paper An integral of Ramanujan and
orthogonal polynomials \jour J. Indian Math. Soc. \vol 51 \yr
1987 \pages 27--36
\endref


\ref\key B1 \by A. Borodin \paper Harmonic analysis on the
infinite symmetric group and the Whittaker kernel \jour
St.~Petersburg Math. J. \vol 12 \yr 2001 \issue 5 \pages 733-759
\endref

\ref\key B2 \by A. Borodin \paper Riemann--Hilbert problem and
the discrete Bessel kernel \jour Intern. Math. Research Notices
\yr 2000 \issue 9 \pages 467--494; {\tt arXiv:\,math.CO/9912093}
\endref


\ref\key BD \by A.~Borodin and P.~Deift \paper Fredholm
determinants, Jimbo--Miwa--Ueno tau--functions, and
representation theory \jour Commun. Pure Appl. Math. \vol 55 \yr
2002 \issue 9 \pages 1160--1230; {\tt arXiv:\, math-ph/0111007}
\endref

\ref\key BO1 \by A.~Borodin and G.~Olshanski \paper Point
processes and the infinite symmetric group. Part III: Fermion
point processes \paperinfo Preprint, 1998, {\tt
arXiv:\,math.RT/9804088}
\endref

\ref\key BO2 \by A.~Borodin and G.~Olshanski \paper Point
processes and the infinite symmetric group \jour Math. Research
Lett. \vol 5 \yr 1998 \pages 799--816; {\tt arXiv:\,
math.RT/9810015}
\endref

\ref\key BO3 \by A.~Borodin and G.~Olshanski \paper
Distributions on partitions, point processes and the
hypergeometric kernel \jour Comm. Math. Phys. \vol 211 \yr 2000
\issue 2 \pages 335--358; {\tt arXiv:\, math.RT/9904010}
\endref

\ref\key BO4 \by A.~Borodin and G.~Olshanski \paper Z--Measures
on partitions, Robinson--Schensted--Knuth correspondence, and
$\beta=2$ random matrix ensembles \inbook In: Random matrix
models and their applications (P.~Bleher and A.~Its, eds).
Cambridge University Press. Mathematical Sciences Research
Institute Publications {\bf 40}, 2001, 71--94; {\tt arXiv:\,
math.CO/9905189}
\endref

\ref\key BO5 \by A.~Borodin and G.~Olshanski \paper Infinite
random matrices and ergodic measures \jour Comm. Math. Phys \vol
223 \yr 2001 \issue 1 \pages 87--123; {\tt arXiv:\,
math-ph/0010015}
\endref

\ref\key BO6 \by A.~Borodin and G.~Olshanski \paper Harmonic
analysis on the infinite--dimensional unitary group and
determinantal point processes \jour Ann. Math., accepted \pages
{\tt arXiv:\, math.RT/0109194}
\endref

\ref\key BO7 \by A.~Borodin and G.~Olshanski \paper Random
partitions and the Gamma kernel \jour Adv. Math. \pages in
press, online publication 2004; {\tt arXiv:\,math-ph/0305043}
\endref

\ref \key Boy \by R.~P.~Boyer \paper Infinite traces of
AF--algebras and characters of $U(\infty)$ \jour J.\ Operator
Theory \vol 9 \yr 1983 \pages 205--236
\endref


\ref\key DVJ \by D.~J.~Daley and D.~Vere--Jones \book An
introduction to the theory of point processes \bookinfo Springer
series in statistics \publ Springer \yr 1988
\endref

\ref\key De \by P.~Deift \paper Integrable operators \inbook In:
Differential operators and spectral theory: M. Sh. Birman's 70th
anniversary collection (V.~Buslaev, M.~Solomyak, D.~Yafaev,
eds.) \bookinfo American Mathematical Society Translations, ser.
2, v. 189 \publ Providence, R.I.: AMS \yr 1999 \pages 69--84
\endref

\ref \key Ed \by A.~Edrei \paper On the generating function of a
doubly--infinite, totally positive sequence \jour Trans.\ Amer.\
Math.\ Soc.\ \vol 74 \issue 3 \pages 367--383 \yr 1953
\endref

\ref\key Fo \by P.~J.~Forrester \book Log--gases and random
matrices \bookinfo Book in preparation, see Forrester's home
page at {\tt http://www.ms.unimelb.edu.au/\~matpjf/matpjf.html}
\endref

\ref\key He \by S.~Helgason \book Groups and geometric analysis.
Integral geometry, invariant differential operators, and
spherical functions \bookinfo Mathematical Surveys and
Monographs {\bf 83} \publ  American Mathematical Society
\publaddr Providence, R.I. \yr 2000
\endref

\ref\key IIKS \by A.~R.~Its, A.~G.~Izergin, V.~E.~Korepin,
N.~A.~Slavnov \paper Differential equations for quantum
correlation functions \jour Intern. J. Mod. Phys. \vol B4 \yr
1990 \pages 10037--1037
\endref

\ref\key JMMS \by M.~Jimbo, T.~Miwa, Y.~M\^ori, and M.~Sato
\paper Density matrix of an impenetrable Bose gas and the fifth
Painlev\'e transcendent \jour Physica D \vol 1 \yr 1980 \pages
80--158
\endref

\ref\key Jo1 \by K.~Johansson \paper Shape fluctuations and
random matrices \jour Commun. Math. Phys. \vol 209 \yr 2000
\issue 2 \pages 437--476; {\tt arXiv:\, math.CO/9903134}
\endref

\ref\key Jo2 \by K.~Johansson \paper Discrete orthogonal
polynomial ensembles and the Plancherel measure \jour Ann. of
Math. (2) \vol 153 \yr 2001 \issue 1 \pages 259--296; {\tt
arXiv:\, math.CO/9906120}
\endref

\ref\key Jo3 \by K.~Johansson \paper Non--intersecting paths,
random tilings and random matrices \jour Probab. Theory Related
Fields \vol 123 \yr 2002 \issue  2 \pages 225--280; {\tt
arXiv:\,math.PR/0011250}
\endref

\ref \key KOV1 \by S.~Kerov, G.~Olshanski, and A.~Vershik \paper
Harmonic analysis on the infinite symmetric group. A deformation
of the regular representation \jour Comptes Rend. Acad. Sci.
Paris, S\'er. I \vol 316 \yr 1993 \pages 773--778
\endref

\ref\key KOV2 \by S.~Kerov, G.~Olshanski, and A.~Vershik \paper
Harmonic analysis on the infinite symmetric group \jour Invent.
Math. \pages in press, online publication 2004; {\tt arXiv:\,
math.RT/0312270}
\endref

\ref\key Len \by A.~Lenard \paper Correlation functions and the
uniqueness of the state in classical statistical mechanics \jour
Comm. Math. Phys \vol 30 \yr 1973 \pages 35--44
\endref

\ref\key Les1 \by P.~A.~Lesky \paper Unendliche und endliche
Orthogonalsysteme von Continuous Hahnpolynomen \jour Results in
Math. \vol 31 \yr 1997 \pages 127--135
\endref

\ref\key Les2 \by P.~A.~Lesky \paper Eine Charakterisierung der
kontinuierlichen und diskreten klassischen Orthogonalpolynome
\paperinfo Preprint 98--12, Mathematisches Institut A,
Universitaet Stuttgart (1998)
\endref

\ref\key Na \by M.~A.~Naimark \book Normed algebras \bookinfo
Translated from the second Russian edition (Moscow, Nauka, 1968)
by Leo F. Boron. Third edition. Wolters--Noordhoff Series of
Monographs and Textbooks on Pure and Applied Mathematics.
Wolters--Noordhoff Publishing, Groningen, 1972
\endref

\ref \key Ner \by Yu.~A.~Neretin \paper Hua type integrals over
unitary groups and over projective limits of unitary groups
\jour Duke Math. J. \vol 114 \issue 2 \pages 239--266 ; {\tt
arXiv:\, math-ph/0010014}
\endref

\ref \key Nes \by N.~I.~Nessonov \paper A complete
classification of the representations of $GL(\infty)$ containing
the identity representation of the unitary subgroup \jour
Mathematics USSR -- Sbornik \vol 58 \yr 1987 \pages 127--147
(translation from Mat. Sb. {\bf 130} (1986), No. 2, 131--150)
\endref

\ref \key OkOl \by A.~Okounkov and G.~Olshanski \paper
Asymptotics of Jack polynomials as the number of variables goes
to infinity  \jour Intern. Math. Research Notices \yr 1998
\issue 13 \pages 641--682
\endref

\ref \key Ol1 \by G.~I.~Ol'shanskii \paper Unitary
representations of infinite--dimensional pairs $(G,K)$ and the
formalism of R.\ Howe \jour Soviet Math. Doklady \vol 27 \issue
2 \yr 1983 \pages 290--294 (translation from Doklady AN SSSR
{\bf 269} (1983), 33--36)
\endref

\ref\key Ol2 \by G.~I.~Ol'shanskii \paper Unitary
representations of the group $SO_0(\infty,\infty)$ as limits of
unitary representations of the groups $SO_0(n,\infty)$ as $n\to
\infty$ \jour  Funct. Anal. Appl. \vol 20 \yr 1986 \issue \pages
292--301
\endref

\ref \key Ol3 \by G.~I.~Ol'shanskii \paper Method of holomorphic
extensions in the theory of unitary representations of
infinite--dimensional classical groups\jour Funct. Anal. Appl.
\vol 22\yr 1988 \issue 4 \pages 273--285
\endref

\ref\key Ol4 \by G.~I.~Ol'shanskii \paper Unitary
representations of $(G,K)$--pairs connected with the infinite
symmetric group S($\infty)$ \jour Leningrad Math. J. \vol
1\issue 4 \yr 1990 \pages 983--1014 (translation from Algebra i
Analiz {\bf 1} (1989), No.4, 178--209)
\endref

\ref \key Ol5 \by G.~I.~Ol'shanskii \paper Unitary
representations of infinite--dimensional pairs $(G,K)$ and the
formalism of R.\ Howe \inbook In: Representation of Lie Groups
and Related Topics \eds A.\ Vershik and D.\ Zhelobenko \bookinfo
Advanced Studies in Contemporary Math. {\bf 7} \publ Gordon and
Breach Science Publishers \publaddr New York etc. \yr 1990
\pages 269--463
\endref

\ref\key Ol6 \by G. Olshanski \paper An introduction to harmonic
analysis on the infinite symmetric group \inbook In: Asymptotic
combinatorics with applications to mathematical physics \ed
A.~M.~ Vershik \bookinfo A European mathematical summer school
held at the Euler Institute, St.~Petersburg, Russia, July 9--20,
2001 \publ Springer Lect. Notes Math. {\bf 1815}, 2003,
127--160; {\tt arXiv:\, math.RT/0311369}
\endref

\ref\key Ol7 \by G. Olshanski \paper The problem of harmonic
analysis on the infinite--dimensional unitary group \jour J.
Funct. Anal. \vol 205 \yr 2003 \issue 2\pages 464--524; {\tt
arXiv:\, math.RT/0109193}
\endref

\ref \key Pi1 \by D.~Pickrell \paper Measures on infinite
dimensional Grassmann manifold \jour J.~Func.\ Anal.\ \vol 70
\yr 1987 \pages 323--356
\endref

\ref \key Pi2 \by D.~Pickrell \paper Separable representations
for automorphism group of infinite symmetric spa\-ces \jour
J.~Func.\ Anal.\ \vol 90 \yr 1990 \pages 1--26
\endref

\ref\key So \by A.~Soshnikov \paper Determinantal random point
fields \jour Russian Math. Surveys \vol 55 \yr 2000 \issue 5
\pages 923--975; {\tt arXiv:\, math.PR/0002099}
\endref

\ref\key VK1 \by A.~M.~Vershik and S.~V.~Kerov \paper Asymptotic
theory of characters of the symmetric group \jour Funct. Anal.
Appl. \vol 15 \yr 1981 \pages 246--255
\endref

\ref \key VK2 \by A.~M.~Vershik and S.~V.~Kerov \paper
Characters and factor representations of the infinite unitary
group \jour Soviet Math.\ Doklady \vol 26 \pages 570--574 \yr
1982
\endref

\ref \key Vo \by D.~Voiculescu \paper Repr\'esentations
factorielles de type {\rm II}${}_1$ de $U(\infty)$ \jour J.\
Math.\ Pures et Appl.\ \vol 55 \pages 1--20 \yr 1976
\endref

\ref\key We \by H.~Weyl \book The classical groups, their
invariants and representations \publ Princeton University Press
\yr 1946
\endref

\ref \key Zhe \by D.~P.~Zhelobenko \book Compact Lie groups and
their representations \publ Nauka, Moscow, 1970 (Russian);
English translation: Transl. Math. Monographs {\bf 40}, Amer.
Math. Soc., Providence, R.I., 1973
\endref

\endRefs
\bye